\documentclass[11pt,a4paper,reqno]{article} 
\usepackage{amsfonts,bm,pifont,alltt,amsthm,amscd,epsfig,amsmath,amssymb,enumerate,url,psfrag,txfonts,stmaryrd}

\newtheorem{theorem}{Theorem}

\newtheorem{proposition}[theorem]{Proposition}

\newtheorem{lemma}[theorem]{Lemma}

\newenvironment{zeproof}{\vskip 2mm\noindent \textbf{\textit{Proof~}}}
                    {\hfill $\blacksquare$ \vskip 2mm \noindent}
\newenvironment{zeproofbis}[1]{\vskip 2mm\noindent \textbf{\textit{Proof (#1)~}}}
                    {\hfill $\blacksquare$ \vskip 2mm \noindent}

\def\card{\mathrm{Card}}

\def\eps{\varepsilon}
\def\II{\mbox{\rm  1\kern-0.20em I}}

\newcommand{\E}{\mathop{\hbox{\rm I\kern-0.17em E}}\nolimits}
\renewcommand{\P}{\mathop{\hbox{\rm I\kern-0.17em P}}\nolimits}

\newcommand{\Z}{\mathbb{Z}}

\newcommand{\N}{\mathbb{N}}
\newcommand{\R}{\mathbb{R}}
\newcommand{\T}{\mathcal{T}}
\newcommand{\TT}{\mathbb{T}}

\renewcommand{\o}{\ensuremath{\omega}}
\renewcommand{\O}{\ensuremath{\Omega}}

\begin{document}

\title{Multiple geodesics with the same direction}

\author{\textsc{David Coupier}\footnote{\texttt{david.coupier@math.univ-lille1.fr}}}

\date{Laboratoire Paul Painlev\'e, UMR 8524}

\maketitle

\noindent
{\bf Abstract:} The directed last-passage percolation (LPP) model with independent exponential times is considered. We complete the study of asymptotic directions of infinite geodesics, started by Ferrari and Pimentel \cite{FP}. In particular, using a recent result of \cite{CH2} and a local modification argument, we prove there is no (random) direction with more than two geodesics with probability $1$.\\

\noindent
{\bf Keywords:} geodesic, last-passage percolation, topological end, random tree.\\

\noindent
{\bf AMS subject classification:} 60K35, 82B43.\\

\section{Introduction}

Commonly, geodesics are known as the generalization of straight lines in Euclidean geometry to curved spaces. For example, the geodesics on a sphere are parts of great circles (arcs). In a random environment, geodesics can be defined as paths optimizing a deterministic procedure. In this paper, the random environment is given by a LPP model in which geodesics are paths with maximal time.\\
\indent
Precisely, let us consider i.i.d. random variables $\omega(z)$, $z\in\N^{2}$, exponentially distributed with parameter $1$. The \textit{last-passage time to $z$} is defined by
$$
G(z) = \max_{\gamma} \sum_{z'\in \gamma} \omega(z')
$$
where the above maximum is taken over all directed paths from the origin to $z$. See Martin \cite{Martin3} for a quite complete survey. Almost surely, the maximum $G(z)$ is reached by only one path: this will be the geodesic to $z$. Thus, the collection of all these geodesics provides a random tree rooted at the origin and spanning all $\N^{2}$. It is called the \textit{geodesic tree} and is denoted by $\T$.\\
\indent
Of course, changing the random environment or the procedure affects the geodesics and the infinite graph they generate. One of the first issues about such an infinite graph concerns its number of topological ends, i.e. the number of infinite self-avoiding geodesics from any fixed vertex. Alexander \cite{A} proved the Minimal Spanning Forest is one-ended. The same statement holds for the Directed Spanning Forest: see Coupier and Tran \cite{CT}.\\
\indent
In first-passage percolation (FPP), results differ according to the nature-- discrete or continuous --of the random environment. Hoffman \cite{Hoffman} proved that for a large class of ergodic FPP processes on the lattice $\Z^{2}$, the number of infinite self-avoiding geodesics is a.s. greater than four. It can even be infinite in some cases. When the lattice $\Z^{2}$ is replaced with the location of an homogeneous Poisson point process on $\R^{2}$, Howard and Newman talk about euclidean FPP models. In \cite{HN}, they got accurate estimates on the fluctuations of geodesics leading to conclusions not yet obtained in the lattice context; quoting them, \textit{a.s. every infinite geodesic has an asymptotic direction and every direction has an infinite geodesic}.\\
Since then, this strategy has been successfully applied to the Radial Spanning Tree by Baccelli and Bordenave \cite{baccellibordenave}, and to the LPP model by Ferrari and Pimentel \cite{FP}. These last results are recalled in Theorem \ref{theorem0}.\\
\indent
In this paper, we finish the study of directions of infinite geodesics in the LPP model. Our results strengthen and complete Ferrari and Pimentel's ones.\\
Indeed, Ferrari and Pimentel \cite{FP}, and Martin \cite{M} had already remarked the almost sure existence of random directions with more than one geodesic. We prove (Part 1. of Theorem \ref{Zetheorem}) the random set of directions with two geodesics is a.s. dense in $[0;\frac{\pi}{2}]$ and countable. Furthermore, with probability $1$, there is no direction with more than two geodesics (Part 2. of Theorem \ref{Zetheorem}). Besides, Ferrari and Pimentel have proved that there is a.s. only one infinite geodesic in each deterministic direction belonging to a subset of $[0;\frac{\pi}{2}]$ of full Lebesgue measure. In Part 3. of Theorem \ref{Zetheorem}, we prove this holds for any direction.\\
In a recent work \cite{CH2}, it has been shown that the probability for the three sites $(0,2)$, $(1,1)$ and $(2,0)$ to be crossed by infinite geodesics is equal to $6-8\log2$. Theorem \ref{nCoex} states such a phenomenon associated to the sites $(0,n),(1,n-1),\ldots,(n,0)$ occurs with positive probability.\\
\indent
Our results are based on a local modification of the geodesic tree. First, remark that any given asymptotic property of the geodesic tree depends on times $\o(z)$, $|z|\leq m$ only through the last-passage times $G(z)$, $|z|=m$ (where $m\in\N$ and $|\cdot|$ denotes the $L^{1}$-norm). Hence, Proposition \ref{prop:EcondB>0} allows to replace the geodesic tree on the set $\{|z|\leq m\}$ with any deterministic directed tree, without changing its structure on $\{|z|\geq m\}$.\\
To our knowledge, no result similar to those of Theorem \ref{Zetheorem} has been established for the models previously mentionned in this section (except for Part 3. whenever an isotropy property holds). Actually, the LPP model takes the advantage to be deeply linked to a particle system, namely the Totally Asymmetric Simple Exclusion Process (TASEP). This coupling, due to Rost \cite{Rost} and described in \cite{CH2}, allows to transfer results about the TASEP to the LPP model. The reader can refer to Theorem 1 of \cite{FP}, Theorems 1 and 2 of \cite{CH2}. Parts 2. and 3. of Theorem \ref{Zetheorem} derive from this coupling too.\\

\indent
The paper is organized as follows. Section \ref{sect:results} contains the definition of the directed LPP model and the results. The local modification argument is detailled in Section \ref{sect:localmodif}. Finally, Theorems \ref{Zetheorem} and \ref{nCoex} are proved in Section \ref{proofs}.

\section{Results}
\label{sect:results}

\begin{figure}[!ht]
\begin{center}
\psfrag{0}{$(0,0)$}
\includegraphics[width=5.5cm,height=4.5cm]{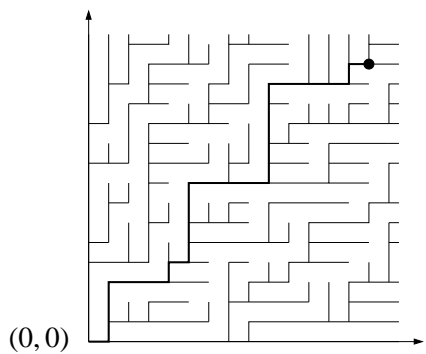}
\end{center}
\caption{\label{fig:geodesictree} An example of the geodesic tree on the set $[0;15]^{2}$. The geodesic of the site $(14,14)$ is represented in bold.}
\end{figure}

Let $\P$ be the Borel probability measure on $\Omega=[0,\infty)^{\N^{2}}$ induced by the family $\{\o(z),z\in\N^{2}\}$ of i.i.d. random variables exponentially distributed with parameter $1$.\\
A \textit{directed path} $\gamma$ from $(0,0)$ to $z$ is a finite sequence of sites $(z_{0},z_{1},\ldots,z_{k})$ with $z_{0}=(0,0)$, $z_{k}=z$ and $z_{i+1}-z_{i}=(1,0)\text{ or } (0,1)$, for $0\leq i\leq k-1$. The quantity
$\sum_{z'\in \gamma} \omega(z')$ represents the time to reach $z$ via $\gamma$. Its maximum over the set $\Gamma(z)$ of all directed paths from $(0,0)$ to $z$ is called the \textit{last-passage time to $z$} and is denoted by $G(z)$:
$$
G(z) = \max_{\gamma\in\Gamma(z)} \sum_{z'\in \gamma} \omega(z') ~.
$$
Since $\P$ is a product measure over $\N^2$ of non-atomic laws, the maximum $G(z)$ is a.s. reached by only one (directed) path, called the \textit{geodesic} to $z$. Let us denote it by $\gamma_{z}$. The collection of all the geodesics $\gamma_{z}$, $z\in\N^{2}$, can be interpreted as a tree rooted at the origin, spanning all the positive quadrant $\N^{2}$ and whose edge set is made up of pairs $\{z',z\}$ with $z-z'=(1,0)\text{ or } (0,1)$ and $z'\in\gamma_{z}$. This random tree is called the \textit{geodesic tree} and is denoted by $\T$. See Figure \ref{fig:geodesictree} for an illustration.

An \textit{infinite geodesic} is merely a infinite branch of the geodesic tree $\T$. Precisely, this is a semi-infinite directed path $(z_{n})_{n\in\N}$ such that, for all $n$, the geodesic $\gamma_{z_{n}}$ is exactly $(z_{0},\ldots,z_{n})$. The existence of infinite geodesics different from the horizontal and vertical axes has been stated by Ferrari and Pimentel. Their results (Propositions 7 and 8 of \cite{FP}) are summarized below. Let us recall that an infinite geodesic $(z_{n})_{n\in\N}$ has direction $\alpha\in[0;\frac{\pi}{2}]$ if
$$
\lim_{n\to\infty} \frac{z_{n}}{\|z_{n}\|} = e^{i \alpha}
$$
(where $\|\cdot\|$ denotes the euclidean norm).

\begin{theorem}[Ferrari and Pimentel \cite{FP}]\hspace*{1cm}
\label{theorem0}
\begin{enumerate}
\item $\P$-a.s. each infinite geodesic has a direction in $[0;\frac{\pi}{2}]$.
\item $\P$-a.s. for all $\alpha\in[0;\frac{\pi}{2}]$, there exists at least one infinite geodesic with direction $\alpha$.
\item There exists a (deterministic) set $D\subset[0;\frac{\pi}{2}]$ of full Lebesgue measure such that, for all $\alpha\in D$, $\P$-a.s. there exists at most one infinite geodesic with direction $\alpha$.
\end{enumerate}
\end{theorem}

The first two parts of Theorem \ref{theorem0} are based on a clever and efficient technic developed by Howard and Newman \cite{HN}. Let $f$ be a positive function on $\R_{+}$. A tree $R$ is said $f$-straight if for all but finitely many vertices $z$ of $R$, the subtree of $R$ coming from $z$ is included in the semi-infinite cone rooted at $z$ with direction $[(0,0);z)$ and angle $f(\|z\|)$. Howard and Newman proved (Proposition 2.8 of \cite{HN}), if $R$ is $f$-straight with $f(x)\to 0$ as $x\to\infty$, then the two assertions corresponding to Parts 1. and 2. of Theorem \ref{theorem0} are satisfied. Thenceforth, an upperbound for the fluctuations of (finite) geodesics (Lemma 10 of \cite{FP}) allows to prove that the geodesic tree $\T$ is a.s. $f$-straight for a suitable function $f$.\\
Besides, Part 3. of Theorem \ref{theorem0} derives from an appropriate use of Fubini's theorem.\\

Our main result completes the study of directions of infinite geodesics of $\T$.

\begin{theorem}\hspace*{1cm}
\label{Zetheorem}
\begin{enumerate}
\item $\P$-a.s. the set of directions with two geodesics is dense in $[0;\frac{\pi}{2}]$ and countable.
\item $\P$-a.s. there is no direction $\alpha\in[0;\frac{\pi}{2}]$ with more than two geodesics.
\item For all $\alpha\in[0;\frac{\pi}{2}]$, $\P$-a.s. there is exactly one geodesic with direction $\alpha$.
\end{enumerate}
\end{theorem}

Let us first remark that multiplying the $\o(z)$'s by a same factor $\lambda$ does not alter the graph structure of $\T$. As a consequence, Theorem \ref{Zetheorem} (and also Theorem \ref{nCoex} below) remains true when replacing the parameter $1$ of the exponential distribution of the $\o(z)$'s with any positive real number.\\
Parts 2. and 3. of Theorem \ref{Zetheorem} rely on the local modification argument of the geodesic tree (see the next section) and on asymptotic results of the LPP model, namely Theorem 1 of \cite{FP} and Theorem 2 of \cite{CH2}.\\
Conversely, the density result of the set of directions with two geodesics is only based on geodesic arguments and Theorem \ref{theorem0}. Let us add this set is necessarily random.\\

This section ends with a coexistence result. Let $n$ be a positive integer. Let $\T_{z}$ be the subtree of $\T$ rooted at $z$; $\T_{z}$ is the collection of geodesics passing by $z$ at which the common part from the origin to $z$ is deleted. We will say there is $n$-\textit{coexistence} if the vertex sets of the $n$ subtrees $\T_{(0,n\!-\!1)},\T_{(1,n\!-\!2)},\ldots,\T_{(n\!-\!1,0)}$ are simultaneously unbounded.\\
When $n=2$, it has been proved (Theorem 1 of \cite{CH2}) that $2$-coexistence occurs with probability $6-8\log2$.

\begin{theorem}
\label{nCoex}
For all $n$, the $n$-coexistence occurs with positive probability.
\end{theorem}

\section{Local modification of the geodesic tree}
\label{sect:localmodif}

Let us denote by $|\cdot|$ the $L^{1}$-norm: $|(x,y)|=x+y$. In this section, we focus on events depending on times $\o(z)$, for $|z|<m$, only through the geodesic tree they generate. Let us start with specifying such events. Let $m\in\N^{\ast}$. A \textit{directed tree} on $\{|z|\leq m\}$ is a graph whose vertex set is $\{|z|\leq m\}$ and whose edge set contains
\begin{itemize}
\item all pairs $\{(x,0),(x+1,0)\}$ and $\{(0,x),(0,x+1)\}$ for $x=0,1,\ldots,m-1$;
\item either $\{z-(1,0),z\}$ or $\{z-(0,1),z\}$, for any $z\in(\N^{\ast})^{2}$ such that $|z|\leq m$.
\end{itemize}
The (finite) set of directed trees on $\{|z|\leq m\}$ is denoted by $\TT_{m}$.\\
Thus, let us denote by $\T_{m}$ the restriction of the geodesic tree $\T$ to the set $\{|z|\leq m\}$: $\T_{m}$ is given by the collection of geodesics $\gamma(z)$, for $|z|\leq m$. With probability $1$, the random tree $\T_{m}$ is $\TT_{m}$-valued. Actually, $\T_{m}$ can take any value of $\TT_{m}$ with positive probability (see Lemma \ref{lem:Tm-dep} at the end of this section). Moreover, for any $T\in\TT_{m}$, the set $\{\o\in\O, \T_{m}(\o)=T\}$ can be expressed as a conjunction of conditions of type ``$G(z)-G(z+(1,-1))$ is positive/negative'' for any $|z|<m$ with positive ordinate. So, this event does not depend on times $\o(z)$ of the diagonal $\{|z|=m\}$.\\
For any nonempty subset $U$ of $\TT_{m}$, the event
$$
\{ \o \in \O , \; \T_{m}(\o) \in U \}
$$
of $\sigma(\o(z), |z|<m)$ is said \textit{$\TT_{m}$-dependent}.\\
Now, we can state the main result of this section. Let us consider an event $S$, with positive probability, describing some asymptotic properties of the geodesic tree $\T$, as for instance ``having more than two infinite geodesics with the same direction''. Then, Proposition \ref{prop:EcondB>0} allows to locally change $\T$-- in order to get additional properties --so that the altered event still occurs with positive probability.

\begin{proposition}
\label{prop:EcondB>0}
Let $m\in\N^{\ast}$ and $S$ be an event of $\sigma( G(z) , \mbox{ for } |z| = m \,;\, \o(z) , \mbox{ for } |z| > m )$ with positive probability. Then, for all $\TT_{m}$-dependent event $B$,
$$
\P ( B \cap S ) > 0 ~.
$$
\end{proposition}

The rest of this section is devoted to the proof of Proposition \ref{prop:EcondB>0}. For that purpose, let us denote by $\mathbf{G}_{m}\in\R_{+}^{m+1}$ the vector of last-passage times on the diagonal $\{|z|=m\}$:
$$
\mathbf{G}_{m} = \left( G(m,0), G(m-1,1),\ldots, G(0,m) \right) ~.
$$
Proposition \ref{prop:EcondB>0} is an immediate consequence of the next result.

\begin{lemma}
\label{lem:AcapB}
Let $m\in\N^{\ast}$ and $B$ be a $\TT_{m}$-dependent event. Then, for any event $A\in\sigma(\mathbf{G}_{m})$, $\P(A)>0$ implies $\P(A\cap B)>0$.
\end{lemma}

\begin{zeproofbis}{of Proposition \ref{prop:EcondB>0}}
Let $m$ be a positive integer. Let $B$ be a $\TT_{m}$-dependent event. Applying Lemma \ref{lem:AcapB} to the event $A=\{\E[\II_{B}\;|\;\mathbf{G}_{m}]=0\}$ which belongs to $\sigma(\mathbf{G}_{m})$ provides
\begin{equation}
\label{>0as}
\E [ \II_{B} \;|\; \mathbf{G}_{m} ] > 0 \; \mbox{ a.s.}
\end{equation}
This means that whatever the value of the vector $\mathbf{G}_{m}$, the tree $\T_{m}$ has a positive probability to satisfy the $\TT_{m}$-dependent event $B$. Moreover, $\P(S)>0$ forces the conditional expectation $\E[\II_{S}\;|\;\mathbf{G}_{m}]$ to be nonzero. Combining with (\ref{>0as}), we get the product $\E[\II_{B}\;|\;\mathbf{G}_{m}] \E[\II_{S}\;|\;\mathbf{G}_{m}]$ is nonzero. So does its expectation:
\begin{equation}
\label{productEspCond}
\E \left\lbrack \E[\II_{B}\;|\;\mathbf{G}_{m}] \E[\II_{S}\;|\;\mathbf{G}_{m}] \right\rbrack > 0 ~.
\end{equation}
To conclude, it remains to adduce the independence of events $B$ and $S$ conditionally to $\mathbf{G}_{m}$. Indeed, conditionally to $\mathbf{G}_{m}$, the event $S$ only depends on times on $\{|z|>m\}$ while $B$ depends on times on $\{|z|<m\}$.
\end{zeproofbis}

Before proving Lemma \ref{lem:AcapB}, let us introduce some notations. Let $\mathbf{X}_{m}\in\R_{+}^{m+1}$ be the vector of times of the diagonal $\{|z|=m\}$
$$
\mathbf{X}_{m} = \left( \o(m,0), \o(m-1,1),\ldots, \o(0,m) \right)
$$
and let $\mathbf{Y}_{m}\in\R_{+}^{m+1}$ such that $\mathbf{G}_{m}=\mathbf{Y}_{m}+\mathbf{X}_{m}$. In the sequel, we will always use bold letters for elements of $\R_{+}^{m+1}$. The random vectors $\mathbf{G}_{m}$, $\mathbf{X}_{m}$ and $\mathbf{Y}_{m}$ respectively induce the probability measures $\P_{\mathbf{G}_{m}}$, $\P_{\mathbf{X}_{m}}$ and $\P_{\mathbf{Y}_{m}}$ on $\R_{+}^{m+1}$, which is endowed with its Borel $\sigma$-algebra $\mathcal{B}_{m+1}$.\\
By hypothesis on the $\o(z)$'s, the random vectors $\mathbf{X}_{m}$ and $\mathbf{Y}_{m}$ are independent. Moreover, $\P_{\mathbf{X}_{m}}$ is absolutely continuous with respect to the Lebesgue measure $\lambda^{m+1}$ on $\R_{+}^{m+1}$. Then, so does for $\P_{\mathbf{G}_{m}}$:
\begin{equation}
\label{abscont}
\lambda^{m+1}(A) =0 \; \mbox{ implies } \; \P_{\mathbf{G}_{m}}(A) = 0 ~.
\end{equation}

\begin{zeproofbis}{of Lemma \ref{lem:AcapB}}
Consider two events $A\in\sigma(\mathbf{G}_{m})$ and $A'\in\mathcal{B}_{m+1}$ such that $A=\mathbf{G}_{m}^{-1}(A')$ and $\P(A)=\P_{\mathbf{G}_{m}}(A')$ is positive. For any $\eps>0$, let $A'_{\eps}$ be the event of $\mathcal{B}_{m+1}$ defined by
$$
A'_{\eps} = A' \cap \left\{ \mathbf{z} \in \R_{+}^{m+1} , \; \inf_{\mathbf{x} \in [\mathbf{z}-\eps;\mathbf{z}]} \mathbf{f}_{m}(\mathbf{x}) \geq \eps \right\} ~,
$$
where $\mathbf{f}_{m}$ denotes the density of $\P_{\mathbf{X}_{m}}$ with respect to $\lambda^{m+1}$, and $[\mathbf{z}-\eps;\mathbf{z}]$ the hyperrectangle
$$
\prod_{i=1}^{m+1} [z_{i} - \eps ; z_{i}] ~,
$$
with $\mathbf{z}=(z_{1},\ldots,z_{m+1})$. Since the density $\mathbf{f}_{m}$ is componentwise decreasing, $A'_{\eps}$ can be written as
$$
A' \cap [\eps ; +\infty[^{m+1} \cap \,\mathbf{f}_{m}^{-1}\left([\eps ; +\infty[\right) ~.
$$
The sequence $(A'_{\eps})_{\eps>0}$ increases as $\eps\searrow 0$. Hence,
$$
\lim_{\eps\to 0} \P_{\mathbf{G}_{m}}(A'_{\eps}) = \P_{\mathbf{G}_{m}} \left( \underset{\eps>0}{\cup} A'_{\eps} \right) = \P_{\mathbf{G}_{m}} \left( A' \cap (\R_{+}^{\ast})^{m+1} \right) = \P_{\mathbf{G}_{m}} \left( A' \right) ~.
$$
About the above equalities, the second one comes from the fact that the density of the exponential distribution is positive at each point of $\R_{+}$. The third one derives from (\ref{abscont}). So, let $\eps>0$ small enough such that $\P_{\mathbf{G}_{m}}(A'_{\eps})$ is positive. A second use of (\ref{abscont}) gives
\begin{equation}
\label{lambdaAeps}
\lambda^{m+1}(A'_{\eps}) > 0 ~.
\end{equation}
Now, let $B$ be a $\TT_{m}$-dependent event. Then we write
\begin{eqnarray*}
\P( A\cap B ) & \geq & \P \left( \mathbf{G}_{m}^{-1}(A'_{\eps}) \cap B \cap \mathbf{Y}_{m}^{-1}([0;\eps]^{m+1}) \right) \\
& = & \int_{\left(\R_{+}^{m+1}\right)^{2}} \E \left\lbrack \II_{B \cap \mathbf{Y}_{m}^{-1}([0;\eps]^{m+1})} \;|\; \mathbf{X}_{m} , \mathbf{Y}_{m} \right\rbrack \, \II_{A'_{\eps}}(\mathbf{x}+\mathbf{y}) \, \mathrm{d}\! \P_{(\mathbf{X}_{m},\mathbf{Y}_{m})} (\mathbf{x},\mathbf{y}) ~.
\end{eqnarray*}
Thanks to the independence between the random vectors $\mathbf{X}_{m}$ and $\mathbf{Y}_{m}$, and the fact that the event $B$ belongs to $\sigma(\o(z), |z|<m)$, the last integral becomes :
$$
\int_{\R_{+}^{m+1}} \left( \int_{\R_{+}^{m+1}} \II_{A'_{\eps}}(\mathbf{x}+\mathbf{y}) \, \mathbf{f}_{m}(\mathbf{x}) \, \mathrm{d} \lambda^{m+1}(\mathbf{x}) \right) \E \left\lbrack \II_{B} \;|\; \mathbf{Y}_{m} \right\rbrack \II_{[0;\eps]^{m+1}}(\mathbf{y}) \, \mathrm{d}\! \P_{\mathbf{Y}_{m}}(\mathbf{y}) ~.
$$
Let $\mathbf{y}\in[0;\eps]^{m+1}$ and $\mathbf{x}\in\R_{+}^{m+1}$. The event $A'_{\eps}$ has been constructed in order to $\mathbf{f}_{m}(\mathbf{x})$ is larger than $\eps$ whenever $\mathbf{x}+\mathbf{y}$ belongs to $A'_{\eps}$. Hence,
$$
\int_{\R_{+}^{m+1}} \II_{A'_{\eps}}(\mathbf{x}+\mathbf{y}) \, \mathbf{f}_{m}(\mathbf{x}) \, \mathrm{d} \lambda^{m+1}(\mathbf{x}) \geq \eps \lambda^{m+1}(A'_{\eps}-\mathbf{y}) = \eps \lambda^{m+1}(A'_{\eps})
$$
since $A'_{\eps}-\mathbf{y}$ is still included in $\R_{+}^{m+1}$. Combining the previous inequalities, we get
\begin{equation}
\label{minorationPAcapB}
\P( A\cap B ) \geq \eps \lambda^{m+1}(A'_{\eps}) \P \left( B \cap \mathbf{Y}_{m}^{-1}([0;\eps]^{m+1}) \right) ~.
\end{equation}
Finally, Lemma \ref{lem:AcapB} follows from (\ref{lambdaAeps}), (\ref{minorationPAcapB}) and Lemma \ref{lem:Tm-dep} stated below.
\end{zeproofbis}

\begin{lemma}
\label{lem:Tm-dep}
Let $B$ be a $\TT_{m}$-dependent event and let $\eps>0$. Then,
\begin{equation}
\label{Bcapeps}
\P \left( B \cap \mathbf{Y}_{m}^{-1}([0;\eps]^{m+1}) \right) > 0 ~.
\end{equation}
\end{lemma}

\begin{zeproof}
It suffices to prove (\ref{Bcapeps}) with $B=\{\o\in\O, \T_{m}(\o)=T\}$ for any $T$ in $\TT_{m}$. Let us proceed by induction. For any integer $1\leq k\leq m$, let $\mathcal{P}[k]$ be the following property:
$$
\exists \eps_{k} > 0 , \; \P \left( \T_{k} = T_{k} \cap \mathbf{Y}_{k}^{-1}([0;\eps_{k}]^{k+1}) \right) > 0 ~,
$$
where $T_{k}$ denotes the restriction of the directed tree $T$ to the set $\{|z|\leq k\}$. Assume $\mathcal{P}[k]$ holds for $k<m$ and let us prove $\mathcal{P}[k+1]$. The edges of $T_{k+1}\setminus T_{k}$ are determined by the signs of differences $G(z)-G(z+(1,-1))$, with $|z|=k$ and $z\not=(k,0)$. So, we have to choose (with positive probability) the vector $\mathbf{X}_{k}$, i.e. the $k+1$ times $\o(0,k),\ldots,\o(k,0)$, satisfying these signs. Here is a way to do it. Pick $\o(0,k)$ in $]2k\eps_{k};(2k+1)\eps_{k}[$. This event, say $C_{0}$, occurs with positive probability. Now:
\begin{itemize}
\item[$\bullet$] If the edge $\{(0,k),(1,k)\}$ belongs to $T_{k+1}\setminus T_{k}$ then pick $\o(1,k-1)$ in $](2k-2)\eps_{k};(2k-1)\eps_{k}[$. Hence,
$$
G(0,k) - G(1,k-1) = \o(0,k) - \o(1,k-1) + Y(0,k) - Y(1,k-1) > 0
$$
since $\o(0,k)-\o(1,k-1)>\eps_{k}$ and $Y(0,k),Y(1,k-1)\in [0;\eps_{k}]$ by hypothesis.
\item[$\bullet$] Otherwise, this is $\{(1,k-1),(1,k)\}$ which belongs to $T_{k+1}\setminus T_{k}$. In this case, pick $\o(1,k-1)$ in $](2k+2)\eps_{k};(2k+3)\eps_{k}[$. This choice ensures that the difference $G(0,k)-G(1,k-1)$ is negative.
\end{itemize}
In both cases, the condition on the time $\o(1,k-1)$ is an event, say $C_{1}$, with positive probability. Thus, we repeat the procedure until the (suitable) choice of $\o(k,0)$. This procedure allows two different times $\o(z)$ and $\o(z')$ of the diagonal $|z|=k$ to belong to the same interval $]j\eps_{k};(j+1)\eps_{k}[$, provided $z$ and $z'$ are not consecutive. When they are, i.e. $z-z'=(1,-1)$, $|\o(z)-\o(z')|$ is larger than $\eps_{k}$ (and smaller than $3\eps_{k}$).\\
The procedure produces $k+1$ events $C_{0},C_{1},\ldots,C_{k}$ with positive probability which are, together with $\{\T_{k}=T_{k}\}\cap\mathbf{Y}_{k}^{-1}([0;\eps_{k}]^{k+1})$, mutually independent. So, by $\mathcal{P}[k]$, the event
$$
\{ \T_{k} = T_{k} \} \cap \mathbf{Y}_{k}^{-1}([0;\eps_{k}]^{k+1}) \cap \left( \overset{k}{\underset{i=0}{\cap}} C_{i} \right)
$$
has a positive probabilty. Moreover, on this event, the times $\o(z)$, $|z|=k$, have been chosen so that the restricted geodesic tree $\T_{k+1}$ coincides with $T_{k+1}$. They also belong to the interval $[0;(4k+1)\eps_{k}]$. Hence, the largest coordinate of the vector $\mathbf{G}_{k}=\mathbf{Y}_{k}+\mathbf{X}_{k}$ is smaller than $(4k+2)\eps_{k}$. So does for $\mathbf{Y}_{k+1}$. Finally, $\mathcal{P}[k+1]$ holds with $\eps_{k+1}=(4k+2)\eps_{k}$.\\
The property $\mathcal{P}[1]$ is true for every $\eps_{1}>0$:
$$
\P \left( \T_{1} = T_{1} \cap \mathbf{Y}_{1}^{-1}([0;\eps_{1}]^{2}) \right) = \P \left( \o(0,0) \leq \eps_{1} \right) > 0 ~.
$$
Hence, the induction starts and gives $\mathcal{P}[m]$ where
$$
\eps_{m} = \eps_{1} \prod_{k=1}^{m-1} (4k+2) < \eps ~,
$$
for $\eps_{1}>0$ small enough.
\end{zeproof}

\section{Proofs}
\label{proofs}

\subsection{Proof of Theorem \ref{Zetheorem}}

Let us start with proving Part 2. of Theorem \ref{Zetheorem}. The idea of the proof can be summarized as follows. From three infinite geodesics the local modification argument allows to assume that they respectively go by sites $(0,2)$, $(1,1)$ and $(2,0)$. The subtree $\T_{(1,1)}$ rooted at $(1,1)$ is then unbounded. Moreover, if these three infinite geodesics have the same direction, $\T_{(1,1)}$ has a null density with respect to the positive quadrant $\N^{2}$. Such a situation is forbidden by Theorem 2 of \cite{CH2}.\\
\indent
Let us proceed by contradiction: assume there exist with positive probability three infinite geodesics with the same (random) direction. Hence, we can find two (deterministic) integers $m>k>0$ such that the event
$$
S(m,k) = \left\{\begin{array}{c}
\mbox{ there exist three infinite geodesics with the same direction such that }\\
\mbox{ only the middle one crosses the diagonal $\{|z|=m\}$ on $(k,m-k)$ }
\end{array}
\right\}
$$
occurs with positive probability. Note that by planarity two geodesics cannot cross each other. So, among the three geodesics mentionned in $S(m,k)$, a middle one can be identified. Moreover, it is crucial to remark that the event $S(m,k)$ depends on times $\o(z)$, for $|z|\leq m$, only through the vector of last-passage times $\mathbf{G}_{m}$. In other words,
$$
S(m,k) \in \sigma( G(z) , \mbox{ for } |z| = m \,;\, \o(z) , \mbox{ for } |z| > m ) ~.
$$
Thus, let us consider the event $B(m,k)$ defined by
$$
B(m,k) = \left\{ \mbox{ $(k,m-k)$ is the only site of $\{|z|=m\}$ whose geodesic goes by $(1,1)$ } \right\} ~.
$$
In other words, $B(m,k)$ means the intersection between the vertex set of $\T_{(1,1)}$ and the diagonal $\{|z|=m\}$ is reduced to $(k,m-k)$. This event only describes the graph structure of $\T_{m}$; it is $\TT_{m}$-dependent (various directed trees of $\TT_{m}$ satisfy $B(m,k)$). Then, Proposition \ref{prop:EcondB>0} applies and gives
\begin{equation}
\label{Esp>0}
\P \left( B(m,k) \cap S(m,k) \right) > 0 ~.
\end{equation}
Let us denote by $C(1,1)$ the vertex set of $\T_{(1,1)}$. On $B(m,k)\cap S(m,k)$, the set $C(1,1)$ is unbounded since it contains all the sites of the middle infinite geodesic of $S(m,k)$. It is also trapped between two geodesics with the same direction. That forces $C(1,1)$ to have a null density with respect to $\N^{2}$:
$$
\lim_{n\to\infty} \frac{1}{n^{2}} \card \left( C(1,1) \cap [0,n]^{2} \right) = 0 ~.
$$
Now, such a situation never happens by Theorem 2 of \cite{CH2}. This contradicts (\ref{Esp>0}).\\

The proof of Part 3. of Theorem \ref{Zetheorem} is also based on the local modification argument. Indeed, Proposition \ref{prop:EcondB>0} allows to consider each direction with two infinite geodesics as a possible direction of the \textit{competition interface} setting out the subtrees $\T_{(1,0)}$ and $\T_{(0,1)}$.\\
\indent
Let $\alpha\in[0;\frac{\pi}{2}]$. Assume that with positive probability there exist two infinite geodesics with direction $\alpha$. So, we can find three integers $0\leq k<k'\leq m$ such that, with positive probability, there exist two infinite geodesics with direction $\alpha$ intersecting the diagonal $\{|z|=m\}$ on $(k,m-k)$ and $(k',m-k')$. Let us denote by $S(m,k,k')$ this event. It belongs to the $\sigma$-algebra generated by the $G(z)$, for $|z|=m$, and the $\o(z)$, for $|z|>m$.\\
Now, let us consider the $\TT_{m}$-dependent event $B(m,k,k')$ for which the geodesics of sites $(k,m-k)$ and $(k',m-k')$ go through respectively $(0,1)$ and $(1,0)$. Using Proposition \ref{prop:EcondB>0}, it follows
\begin{equation}
\label{Proba>0}
\P \left( B(m,k,k') \cap S(m,k,k') \right) > 0 ~.
\end{equation}
Let us respectively denote by $C(0,1)$ and $C(1,0)$ the vertex sets of the subtrees $\T_{(0,1)}$ and $\T_{(1,0)}$. These two sets provide a random partition of $\N^{2}\setminus\{(0,0)\}$. In \cite{FP}, Ferrari and Pimentel studied the asymptotic behaviour of the boundary between $C(0,1)$ and $C(1,0)$. This boundary is modeled by an infinite directed path, called the competition interface. Ferrari and Pimentel proved (Theorem 1 of \cite{FP}) that the competition interface has $\P$-a.s. an asymptotic direction in $[0;\frac{\pi}{2}]$, whose distribution has no atom. Now, a contradiction with (\ref{Proba>0}) appears. Indeed, on the event $B(m,k,k')\cap S(m,k,k')$, the competition interface has direction $\alpha$.\\
In conclusion, there is no more than one infinite geodesic with direction $\alpha$ with probability $1$. Part 2. of Theorem \ref{theorem0} completes the proof.\\

It remains to show the first part of Theorem \ref{Zetheorem}, i.e. $\P$-a.s. the set of directions with two geodesics is dense in $[0;\frac{\pi}{2}]$ and countable. Let $I$ be a nonempty interval of $[0;\frac{\pi}{2}]$ and let $\alpha,\beta\in I$ such that $\alpha<\beta$. By Part 2. of Theorem \ref{theorem0}, there exists with probability $1$ two infinite geodesics with direction $\alpha$ and $\beta$, say respectively $\gamma_{\alpha}$ and $\gamma_{\beta}$. These two infinite paths have a common part, say $z_{0}=(0,0),z_{1},\ldots,z_{k}$ (possibly reduced to the origin if $z_{k}=z_{0}$). Thus, they bifurcate at $z_{k}$: $z_{k}+(1,0)$ belongs to $\gamma_{\alpha}$, and $z_{k}+(0,1)$ to $\gamma_{\beta}$.\\
Let us define the highest infinite path passing by $z_{k+1}:=z_{k}+(1,0)$. This path denoted by $\gamma_{\alpha}^{+}$ is inductively built as follows. If $z_{n}\in\gamma_{\alpha}^{+}$, $n\geq k+1$, has exactly one child, say $z_{n}'$, with infinitely many descendant in the geodesic tree, then we put $z_{n+1}=z_{n}'$. Otherwise (it has two such children), we put $z_{n+1}=z_{n}+(0,1)$. In the same way, we define the lowest infinite path passing by $z_{k}+(0,1)$; $\gamma_{\beta}^{-}$. Part 1. of Theorem \ref{theorem0} says that $\gamma_{\alpha}^{+}$ and $\gamma_{\beta}^{-}$ almost surely have asymptotic directions. Say respectively $\alpha^{+}$ and $\beta^{-}$. By construction and planarity,
$$
\alpha \leq \alpha^{+} \leq \beta^{-} \leq \beta ~.
$$
Furthermore, the previous construction forces all the paths trapped between $\gamma_{\alpha}^{+}$ and $\gamma_{\beta}^{-}$, and different from those two paths, to be finite. So $\alpha^{+}=\beta^{-}$ ($\in I$). Otherwise, this would contradict Part 2. of Theorem \ref{theorem0}.\\
Hence, the interval $I$ contains with probability $1$ a direction with two infinite geodesics. To conclude, we use the countability and the density of rational numbers into $[0;\frac{\pi}{2}]$: a.s. the set of directions with two infinite geodesics meets every nonempty interval of $[0;\frac{\pi}{2}]$.\\
Finally, let us state the countable character of the set of directions with two infinite geodesics. To do it we consider the random application that associates to each couple of infinite geodesics (the first one above the second one) with the same direction, the site (of $\N^{2}$) at which they bifurcate. Thus, it suffices to remark that this application is $\P$-a.s injective since there is no more than two infinite geodesics with the same direction.

\subsection{Proof of Theorem \ref{nCoex}}

Let $n$ be a positive integer. Consider $n$ different angles $\theta_{1}<\ldots<\theta_{n}$ in $[0;\frac{\pi}{2}]$. We denote by $\gamma_{1},\ldots,\gamma_{n}$ the (different) infinite geodesics with direction respectively $\theta_{1},\ldots,\theta_{n}$. Their existence is almost surely ensured by the third part of Theorem \ref{Zetheorem}. Thus, there exists an integer $m$ such that, with positive probability, the $\gamma_{i}$'s intersect the diagonal $\{|z|=m\}$ on $n$ different sites. In the same way, one can find $n$ sites $z_{1},\ldots,z_{n}$ in $\{|z|=m\}$ such that, with positive probability, for any $1\leq i\leq n$, $\gamma_{i}$ goes by $z_{i}$. Let us denote by $S(m;z_{1},\ldots,z_{n})$ this event.\\
Thus remark that each site $z_{i}$ belongs to the quadrant $(i-1,n-i)+\N^{2}$. So, $z_{i}$ can be reached by a directed path coming from $(i-1,n-i)$. Now, let us consider the event $B(m;z_{1},\ldots,z_{n})$ for which, for any $1\leq i\leq n$, the geodesic to $z_{i}$ goes through $(i-1,n-i)$. This is a $\TT_{m}$-dependent event.\\
On the one hand, $B(m;z_{1},\ldots,z_{n})$ and $S(m;z_{1},\ldots,z_{n})$ satisfy the hypotheses of Proposition \ref{prop:EcondB>0}. Then,
$$
\P \left( B(m;z_{1},\ldots,z_{n}) \cap S(m;z_{1},\ldots,z_{n}) \right) > 0 ~.
$$
On the other hand, on $B(m;z_{1},\ldots,z_{n})\cap S(m;z_{1},\ldots,z_{n})$, an infinite geodesic emanates from each site of the diagonal $\{|z|=n-1\}$. This proves that $n$-coexistence occurs with positive probability.

\section*{Acknowledgments}

The author thanks J.~B. Martin for enlightening discussions and J. Kahn for its useful advice.

\bibliographystyle{plain}
\bibliography{LPP_bibli}

\begin{thebibliography}{10}

\bibitem{A}
K.~S. Alexander.
\newblock Percolation and minimal spanning forest in infinite graphs.
\newblock {\em The Annals of Probability}, 23(1):87--104, 1995.

\bibitem{baccellibordenave}
F.~Bacelli and C.~Bordenave.
\newblock The radial spanning tree of a {P}oisson point process.
\newblock {\em Annals of Applied Probability}, 17(1):305--359, 2007.

\bibitem{CH2}
D.~Coupier and P.~Heinrich.
\newblock Coexistence probability in the last passage percolation model is
  $6-8\log2$.
\newblock arXiv:1007.0652, 2010.

\bibitem{CT}
D.~Coupier and V.~C. Tran.
\newblock The directed spanning forest is almost surely a tree.
\newblock arXiv:1010.0773v2, 2010.

\bibitem{FP}
P.~A. Ferrari and L.~P.~R. Pimentel.
\newblock Competition interfaces and second class particles.
\newblock {\em Ann. Probab.}, 33(4):1235--1254, 2005.

\bibitem{Hoffman}
C.~Hoffman.
\newblock Geodesics in first passage percolation.
\newblock {\em Ann. Appl. Probab.}, 18(5):1944--1969, 2008.

\bibitem{HN}
C.~D. Howard and C.~M. Newman.
\newblock Geodesics and spanning trees for {E}uclidean first-passage
  percolation.
\newblock {\em Annals of Probability}, 29:577--623, 2001.

\bibitem{Martin3}
J.~B. Martin.
\newblock Last-passage percolation with general weight distribution.
\newblock {\em Markov Process. Related Fields}, 12(2):273--299, 2006.

\bibitem{M}
J.~B. Martin.
\newblock Private communication.
\newblock 2011.

\bibitem{Rost}
H.~Rost.
\newblock Nonequilibrium behaviour of a many particle process: density profile
  and local equilibria.
\newblock {\em Z. Wahrsch. Verw. Gebiete}, 58(1):41--53, 1981.

\end{thebibliography}

\end{document}